\documentclass[12pt]{article}
\usepackage{here,amsfonts,amssymb}
\usepackage{amsmath}
\newtheorem{thm}{Theorem}[section]

\newtheorem{prop}[thm]{Proposition}

\newbox\sample
\newif\ifproofmode
\newif\ifsymindex
\global\symindexfalse
\newwrite\inx
\def\indsyma#1#2{\ifproofmode\marginpar{$\scriptstyle#1$}\fi%
\ifx#2\empty\write\inx{$\noexpand#1$,\space\thepage}%
\write\inx{\string\newline}\else%
\write\inx{$\noexpand#1$,\space#2,\space\thepage}%
\write\inx{\string\newline}\fi\ignorespaces}%
\def\indsym#1#2{\ifsymindex%
\ifproofmode\marginpar{$\scriptstyle#1$}\fi%
\ifx#2\empty\write\inx{\string\item \space$\noexpand#1$,\space\thepage}%
\else%
\write\inx{\string\item \space$\noexpand#1$,\space#2,\space\thepage}%
\fi\ignorespaces\fi}%
\newskip\dangerskipb
\newskip\dangerskip
\def\hang{\hangindent\dangerskip}

\def\proof{\noindent{\it Proof\/}.\enspace}

\font\manual=manfnt at 12pt
\def\danbend{{\manual\char127}}
\def\ddatanger{\medbreak\begingroup\clubpenalty=10000
 \def\par{\endgraf\endgroup\medbreak} \noindent\hang\hangafter=-2
 \hbox to0pt{\hskip-3.5pc\danbend\kern1pt%
\danbend\hfill}}

\def\dobdownarrow{\mathop{\vbox{\kern2pt \hbox{$\Big\downarrow$}\kern-16.5pt
                          \nointerlineskip\hbox{$\Big\downarrow$}}}}

\def\lrightarrow{\hbox to 25pt{\rightarrowfill}}

\def\supexp{exp(m,n,p)=m^{m^{m^{\cdot^{\cdot^{\cdot^{m^{p}}}}}}}
\vbox{\hbox{$\Big\}\scriptstyle n$}\kern0pt}}

\def\supexpo#1#2#3{#1^{#1^{\cdot^{\cdot^{\cdot^{#1^{#2}}}}}}
\vbox{\hbox{$\Big\}\scriptstyle #3$}\kern0pt}}

\def\sqr#1#2{{\vcenter{\hrule height .#2pt
         \hbox{\vrule width.#2pt height#1pt \kern#1pt
             \vrule width.#2pt}
         \hrule height.#2pt}}}

\def\bigsquare{\mathchoice\sqr76\sqr76\sqr{2.1}3\sqr{1.5}3}

\newskip\bogcentering \bogcentering= 0pt plus 1000pt minus 1000pt 

\def\matth{\mathsurround=0pt}
\def\fakrightarrowfill{$\matth \mathord- \mkern-6mu
  \cleaders\hbox{$\mkern-2mu \mathord- \mkern-2mu$}\hfill
 \mkern-6mu \mathord\rightarrow$}

\def\fakoverrightarrow#1{\vbox{\ialign{##\crcr
  \fakrightarrowfill\crcr\noalign{\kern-1pt\nointerlineskip}
 $\hfil\displaystyle{#1}\hfil$\crcr}}}

\newif\ifdtatp

\def\displaty{%
\global \dtatptrue \openup \jot \matth \everycr{\noalign{\ifdtatp \global 
\dtatpfalse \vskip -\lineskiplimit \vskip \normallineskiplimit \else 
\penalty \interdisplaylinepenalty \fi }}}

\def\displaylignes#1{\displaty
   \halign{\hbox to\displaywidth{$\displaystyle##$}\crcr
   #1\crcr}}
\def\leqaligneno#1{\displaty \tabskip=\bogcentering
 \halign to\displaywidth{\hfil$\displaystyle{##}$\tabskip=0pt
 &$\displaystyle{{}##}$\hfil\tabskip=\bogcentering
 &\kern-\displaywidth\rlap{$##$}\tabskip=\displaywidthpt\crcr
 #1\crcr}}
\def\ligne{\hbox to\hsize}
\newdimen\nouvpagewidth
\newdimen\offwidth
\newdimen\lawidthoui
\nouvpagewidth=\lawidthoui
\def\kboxit#1{\vbox{\hrule\hbox{\vrule\kern3pt
              \vbox{\kern3pt#1\kern3pt}\kern3pt\vrule}\hrule}}

\def\kboxitb#1{\vbox{\hrule\hbox{\vrule\kern3pt
              \vbox{\kern3pt#1\kern3pt}\kern3pt\vrule}\hrule}}

\def\laboxaround#1{
\aboxaround{\hbox to\hsize{\hfill\box2\hfill}}{#1}
}

\def\boxar#1#2{
\aboxaround{\hbox to\hsize{\hfill#1\hfill}}{#2}
}

\def\aboxaround#1#2{
\setbox4=\vbox{\hsize #2\noindent\strut#1\strut}
\kboxitb{\box4}}

\def\kframeit#1{\vbox{\hrule\hbox{\vrule\kern5pt
              \vbox{\kern5pt#1\kern5pt}\kern5pt\vrule}\hrule}}

\newskip\savnormalbaselineskip
\newskip\savnormallineskip
\newdimen\savnormallineskiplimit

\def\transpos#1{#1^{\top}}
\def\amsdiagmat#1#2{%
\begin{pmatrix}
#1_{1}&        &\ldots& \\
                & #1_{2}  &\ldots& \\
\vdots&\vdots&\ddots&\vdots\\
      &      &\ldots& #1_{#2}
\end{pmatrix}
}
\def\amsmata#1#2#3#4{%
\begin{pmatrix}
#1 & #2\\
#3 & #4
\end{pmatrix}
}

\begin{document}
\title{
Remarks on the Cayley Representation of Orthogonal Matrices
and on Perturbing the Diagonal of a Matrix to Make it Invertible
\\ }
\author{Jean Gallier\\
%
 \\
Department of Computer and Information Science\\
University of Pennsylvania\\
Philadelphia, PA 19104, USA\\
{\tt jean@cis.upenn.edu}
}
\maketitle
\vspace{0.3cm}
\noindent
{\bf Abstract.}
This note contains two remarks. The first remark concerns
the extension of the well-known Cayley representation of rotation
matrices by skew-symmetric matrices
to rotation matrices admitting $-1$ as an eigenvalue
and then to all orthogonal matrices. We review a method due to
Hermann Weyl and another method involving multiplication
by a diagonal matrix whose entries are $+1$ or $-1$.
The second remark has to do with ways of flipping the signs
of the entries of a diagonal matrix, $C$, with nonzero diagonal
entries, obtaining a new matrix, $E$, so that $E + A$
is invertible, where $A$ is any given matrix 
(invertible or not).

\vfill\eject
\section{The Cayley Representation of Orthogonal Matrices}
\label{sec1}
Given any rotation matrix, $R\in \mathbf{SO}(n)$, 
if $R$ does not admit $-1$ as an eigenvalue, then
there is a unique skew-symmetric matrix, $S$,
($\transpos{S} = -S$) so that
\[
R = (I - S)(I + S)^{-1}.
\]
This is a classical result of Cayley \cite{Cayley} (1846)
and $R$ is called the {\it Cayley transform of $S$\/}. 
Among other sources, a proof can be found in
Hermann Weyl's beautiful book {\sl The Classical Groups\/}
\cite{Weyl46}, Chapter II, Section 10, Theorem 2.10.B (page 57).

\medskip
As we can see, this representation misses rotation matrices
admitting the eigenvalue $-1$, and of course, 
as $\det((I - S)(I + S)^{-1}) = +1$, it misses
improper orthogonal matrices, i.e., those matrices
$R\in \mathbf{O}(n)$ with $\det(R) = -1$.

\medskip\noindent
{\bf Question 1}.
Is there a way to extend the Cayley representation to all 
rotation matrices (matrices in $\mathbf{SO}(n)$)?

\medskip\noindent
{\bf Question 2}.
Is there a way to extend the Cayley representation to all 
orthogonal matrices (matrices in $\mathbf{O}(n)$)?

\medskip\noindent
{\bf Answer}: Yes in both cases!

\medskip
An answer to Question 1 is given in Weyl's book \cite{Weyl46},
Chapter II, Section 10, Lemma 2.10.D (page 60):

\begin{prop} (Weyl)
\label{Weyl1}
Every rotation matrix, $R\in \mathbf{SO}(n)$, can be expressed as
a product
\[
R = (I - S_1)(I + S_1)^{-1}(I - S_2)(I + S_2)^{-1},
\]
where $S_1$ and $S_2$ are skew-symmetric matrices.
\end{prop}

\medskip
Thus, if we allow two Cayley representation matrices, we can capture
orthogonal matrices having an even number of $-1$ as eigenvalues.
Actually, proposition \ref{Weyl1} can be sharpened slightly as follows:

\begin{prop} 
\label{prop2}
Every rotation matrix, $R\in \mathbf{SO}(n)$, can be expressed as
\[
R = \Bigl((I - S)(I + S)^{-1}\Bigr)^2
\]
where $S$ is a  skew-symmetric matrix.
\end{prop}

\medskip
Proposition \ref{prop2} can be easily proved using the following
well-known normal form for orthogonal matrices:

\begin{prop} 
\label{prop3}
For every orthogonal matrix, $R\in  \mathbf{O}(n)$,
there is an orthogonal matrix $P$
and a block diagonal matrix $D$
such that $R = PD\,\transpos{P}$, 
where $D$ is of the form
\[
D = \amsdiagmat{D}{p}
\]
such that each block $D_i$ is either $1$, $-1$,
or a two-dimensional matrix of the form
$$D_i = \amsmata{\cos\theta_i}{-\sin\theta_i}{\sin\theta_i}{\cos\theta_i}$$
where $0 < \theta_i < \pi$.
\end{prop}

\medskip
In particular, if $R$ is a rotation matrix ($R\in  \mathbf{SO}(n)$), then
it has an even number of eigenvalues $-1$. So, they can be grouped
into two-dimensional rotation matrices of the form
\[
\amsmata{-1}{0}{0}{-1},
\]
i.e., we allow $\theta_i = \pi$ and we may assume that
$D$ does not contain one-dimensional blocks of the form $-1$.

\medskip
A proof of Proposition \ref{prop3} can be found in
Gantmacher \cite{Gantmacher1}, Chapter IX, Section 13 (page 285),
or Berger \cite{Berger90}, or Gallier \cite{Gallbook2},
Chapter 11, Section 11.4 (Theorem 11.4.5).

\medskip
Now, for every two-dimensional rotation matrix
\[
T = \amsmata{\cos\theta}{-\sin\theta}{\sin\theta}{\cos\theta}
\]
with $0 < \theta \leq \pi$, observe that
\[
T^{\frac{1}{2}} = 
 \amsmata{\cos(\theta/2)}{-\sin(\theta/2)}{\sin(\theta/2)}{\cos(\theta/2)}
\]
does not admit $-1$ as an eigenvalue (since $0 < \theta/2 \leq\pi/2$)
and $T = \left(T^{\frac{1}{2}}\right)^2$. Thus,
if we form the matrix $R^{\frac{1}{2}}$ by replacing each
two-dimensional block $D_i$ in the above normal form
by $D_i^{\frac{1}{2}}$, we obtain a rotation matrix that
does not admit $-1$ as an eigenvalue, 
$R = \left(R^{\frac{1}{2}}\right)^2$ and the Cayley transform of
$R^{\frac{1}{2}}$ is well defined. Therefore, we have proved 
Proposition \ref{prop2}.
$\bigsquare$

\bigskip
Next, why is the answer to Question 2 also yes?

\medskip
This is because
\begin{prop}
\label{prop4}
For any orthogonal matrix, $R\in \mathbf{O}(n)$, there is
some diagonal matrix, $E$,  whose entries are $+1$ or $-1$, 
and some skew-symmetric matrix, $S$, so that
\[
R = E(I - S)(I + S)^{-1}.
\]
\end{prop}

\medskip
As such matrices $E$ are orthogonal, all matrices $E(I - S)(I + S)^{-1}$
are orthogonal, so we have a Cayley-style representation of all 
orthogonal matrices.

\medskip
I am not sure when Proposition \ref{prop4} was discovered 
and originally published.
Since I could not locate this result in Weyl's book \cite{Weyl46},
I assume that it was not known before 1946, but I did stumble on
it as an exercise in Richard Bellman's classic \cite{Bellman}, first
published in 1960,
 Chapter 6, Section 4, Exercise 11, page 91-92 
(see also, Exercises, 7, 8, 9, and 10).

\medskip
Why does this work?

\medskip\noindent
{\bf Fact E}: Because, for every $n\times n$ matrix, $A$ (invertible 
or not), there some diagonal matrix, $E$,  whose entries are $+1$ or $-1$, 
so that $I + EA$ is invertible!

\medskip
This is Exercise 10 in Bellman \cite{Bellman} (Chapter 6, Section 4, page 91).
Using Fact E, it is easy to prove Proposition \ref{prop4}.

\medskip\noindent
{\it Proof of Proposition \ref{prop4}\/}.
Let $R\in \mathbf{O}(n)$ be any orthogonal matrix.
By Fact E, we can find a diagonal matrix, $E$ (with diagonal entries
$\pm1$), so that $I + ER$ is invertible. But then, as $E$ is orthogonal,
$ER$ is an orthogonal matrix that does not admit the eigenvalue $-1$ and
so, by the Cayley representation theorem, there is a
skew-symmetric matrix, $S$, so that
\[
ER = (I - S)(I + S)^{-1}.
\]
However, notice that $E^2 = I$, so we get
\[
R = E(I - S)(I + S)^{-1},
\]
as claimed.
$\bigsquare$

\medskip
But Why does Fact E hold?

\medskip
As we just observed, $E^2 = I$, so by multiplying by $E$,
\[
\hbox{$I + EA$ is invertible iff $E + A$ is.}
\]

\medskip
Thus, we are naturally led to the following problem:
If $A$ is any $n\times n$ matrix, is there a way to perturb
the diagonal entries of $A$, i.e., to add some diagonal
matrix,  $C = \mathrm{diag}(c_1, \ldots, c_n)$, to $A$ so that
$C + A$ becomes invertible?

\medskip
Indeed this can be done, and we will show in the next section
that what matters is not the magnitude of the perturbation
but the signs of the entries being added.

\section{Perturbing the Diagonal of a Matrix to Make it Invertible}
\label{sec2}
In this section we prove the following result:

\begin{prop}
\label{prop5}
For every $n\times n$ matrix (invertible or not), $A$,
and every any diagonal matrix,  $C = \mathrm{diag}(c_1, \ldots, c_n)$,
with $c_i \not= 0$ for $i = 1, \ldots, n$, 
there an assignment of signs, $\epsilon_i = \pm 1$, so
that if $E = \mathrm{diag}(\epsilon_1 c_1, \ldots, \epsilon_n c_n)$, then
$E + A$ is invertible.
\end{prop}

\proof
Let us evaluate the determinant of $C + A$. We see that
$\Delta = \det(C + A)$ is a polynomial of degree $n$ in the variables
$c_1, \ldots, c_n$ and that all the  monomials of $\Delta$ consist
of  products of distinct variables (i.e., every variable occurring in 
a monomial has degree $1$). In particular, $\Delta$ contains the 
monomial $c_1 \cdots c_n$.
In order to prove Proposition \ref{prop5}, it will suffice to prove

\begin{prop}
\label{prop6}
Given any polyomial, $P(x_1, \ldots, x_n)$, of degree $n$ (in the
indeterminates $x_1, \ldots, x_n$  and over any integral domain
of characteristic unequal to $2$), if
every monomial in $P$ is a product of distinct variables,
for every $n$-tuple $(c_1, \ldots, c_n)$ such that
$c_i \not= 0$ for $i = 1, \ldots, n$, then
there is an assignment of signs, $\epsilon_i = \pm 1$, so that
\[
P(\epsilon_1 c_1, \ldots, \epsilon_n c_n) \not= 0.
\]
\end{prop}

\medskip
Clearly, any assignment of signs given by Proposition \ref{prop6}
will make $\det(E + A) \not= 0$, proving Proposition   \ref{prop5}.
$\bigsquare$

\medskip
It remains to prove Proposition  \ref{prop6}.

\bigskip\noindent
{\it Proof of Proposition \ref{prop6}\/}.
We proceed by induction on $n$ (starting with $n = 1$).
For $n = 1$, the polynomial $P(x_1)$ is of the form
$P(x_1) = a + bx_1$, with $b\not= 0$ since $\mathrm{deg}(P) = 1$.
Obviously, for any $c\not= 0$, either $a + bc\not= 0$ or
$a - bc \not= 0$ (otherwise, $2bc = 0$, contradicting
$b\not= 0$, $c\not= 0$ and the ring being an integral domain
of characteristic $\not= 2$).

\medskip
Assume the induction hypothesis holds for any $n\geq 1$ and
let $P(x_1, \ldots, x_{n+1})$ be a polynomial of degree $n+1$
satisfying the conditions of Proposition \ref{prop6}. 
Then, $P$ must be of the form
\[
P(x_1, \ldots, x_n, x_{n+1})  
= Q(x_1, \ldots, x_n) + S(x_1, \ldots, x_n)x_{n+1},
\]
where both $Q(x_1, \ldots, x_n)$  and $S(x_1, \ldots, x_n)$
are polynomials in $x_1, \ldots, x_n$ and  $S(x_1, \ldots, x_n)$
is of degree $n$ and all monomials in
$S$ are products of distinct variables. By the induction
hypothesis, we can find $(\epsilon_1, \ldots, \epsilon_n)$,
with $\epsilon_i = \pm 1$, so that
\[
S(\epsilon_1 c_1, \ldots, \epsilon_n c_n) \not= 0.
\]
But now, we are back to the case $n = 1$ with 
the polynomial
\[
Q(\epsilon_1 c_1, \ldots, \epsilon_n c_n) + 
S(\epsilon_1 c_1, \ldots, \epsilon_n c_n)x_{n+1},
\]
and we can find $\epsilon_{n+1} = \pm 1$ so that
\[
P(\epsilon_1 c_1, \ldots, \epsilon_n c_n, \epsilon_{n+1}c_{n+1}) =
Q(\epsilon_1 c_1, \ldots, \epsilon_n c_n) + 
S(\epsilon_1 c_1, \ldots, \epsilon_n c_n)\epsilon_{n+1}c_{n+1}
\not= 0,
\]
establishing the induction hypothesis.
$\bigsquare$

\medskip
Note that in Proposition \ref{prop5}, the $c_i$ can be made
arbitrarily small or large, as long as they are not zero.
Thus, we see as a corollary that any matrix can be made
invertible by a very small perturbation of its diagonal elements.
What matters is the signs that are assigned to the perturbation.

\medskip
Another nice proof of Fact E is given in a short note by William Kahan
\cite{Kahan2004}. Due to its elegance, we feel compelled to
sketch Kahan's proof.
This proof uses two facts:

\begin{enumerate}
\item[(1)]
If $A = (A_1, \ldots, A_{n-1}, U)$ and $B = (A_1, \ldots, A_{n-1}, V)$
are two $n\times n$ matrices that differ in their last column,
then
\[
\det(A + B) = 2^{n-1}(\det(A) + \det(B)).
\]
This is because determinants are multilinear (alternating) maps of
their columns. Therefore, if $\det(A) = \det(B) = 0$, then
$\det(A + B) = 0$. Obviously, this fact also holds whenever
$A$ and $B$ differ by just one column (not just the last one).
\item[(2)]
For every $k = 0, \ldots 2^{n} -1$, if we write $k$ in binary
as $k = k_n \cdots k_1$, then let $E_k$ be the diagonal matrix
whose $i$th diagonal entry is $-1$ iff $k_i = 1$, else $+1$
iff $k_i = 0$. For example, $E_0 = I$ and $E_{2^n -1} = -I$.
Observe that $E_k$ and $E_{k+1}$ differ by exactly one column.
Then, it is easy to see that
\[
E_0 + E_1 + \cdots + E_{2^n -1} = 0.
\]
\end{enumerate}

The proof proceeds by contradiction. Assume that
$\det(I + E_kA) = 0$, \\
for $k = 0, \ldots, 2^n -1$.
The crux of the proof is that
\[
\det(I + E_0A + I + E_1A + I + E_2A + \cdots + I + E_{2^n -1}A) = 0.
\]
However, as $E_0 + E_1 + \cdots + E_{2^n -1} = 0$, we see that
\[
I + E_0A + I + E_1A + I + E_2A + \cdots + I + E_{2^n -1}A 
= 2^n I,
\]
and so,
\[
0 = \det(I + E_0A + I + E_1A + I + E_2A + \cdots + I + E_{2^n -1}A) =
\det(2^n I) = 2^n \not= 0,
\]
a contradiction!

\medskip
To prove that
$\det(I + E_0A + I + E_1A + I + E_2A + \cdots + I + E_{2^n -1}A) = 0$, we 
observe using fact (2) that,
\[
\det(I + E_{2i}A + I + E_{2i+1}A) = 
\det(I + E_{2i}A) + \det(I + E_{2i+1}A) = 0,
\]
for $i = 0, \ldots, 2^{n - 1} -1$; similarly,
\[
\det(I + E_{4i}A + I + E_{4i+1}A + I + E_{4i+2}A + I + E_{4i+3}A) = 0,
\] 
for $i = 0, \ldots, 2^{n - 2} -1$; by induction, we get
\[
\det(I + E_0A + I + E_1A + I + E_2A + \cdots + I + E_{2^n -1}A) = 0,
\]
which concludes the proof.

\bigskip\noindent
{\bf Final Questions\/}:
\begin{enumerate}
\item[(1)]
When was Fact E first stated and by whom (similarly for
Proposition \ref{prop4})?
\item[(2)]
Can Proposition \ref{prop5} be generalized 
to non-diagonal matrices (in an interesting way)?
\end{enumerate}


\bibliographystyle{plain} 
\end{document}